\documentclass[a4paper,12pt,draft]{article}
\usepackage{amsmath,amsfonts,amssymb,amsthm}
\usepackage{t1enc}

\newcommand{\ex}{{\mathbf{E}}}
\newcommand{\gen}{{\Delta^{\frac{\alpha}{2}}}}
\newcommand{\ind}{{\mathbf{1}}}
\newcommand{\pr}{{\mathbf{P}}}
\newcommand{\R}{{\mathbf{R}}}
\newcommand{\Rd}{{\R^d}}
\newcommand{\sub}{\subseteq}

\DeclareMathOperator{\diam}{diam}

\newtheorem{theorem}{Theorem}
\newtheorem*{lemma}{Lemma}
\newtheorem*{proposition}{Proposition}

\begin{document}
\sloppy


\title{Spectral gap estimate for fractional Laplacian}
\author{Mateusz Kwa\'snicki\thanks{Work supported by KBN grant 1 P03A 026 29} \\ {\small Institute of Mathematics and Computer Science} \\ {\small Wroc{\l}aw University of Technology} \\ {\small Wybrze\.ze Wyspia\'nskiego 27, 50-370 Wroc{\l}aw, Poland}}
\maketitle


\begin{abstract}
A lower bound estimate $\lambda_2 - \lambda_1 \ge c \, (\lambda_1)^{-\frac{d}{\alpha}} \, (\diam D)^{-d - \alpha}$ for the spectral gap of the Dirichlet fractional Laplacian $\gen$ on arbitrary bounded domain $D \sub \Rd$ is proved. This follows from a variational formula for the spectral gap (\cite{bib:dk}) and an upper bound estimate for the supremum norm of the ground state eigenfunction.
\end{abstract}


\section*{Introduction and main results}

Isotropic stable processes and fractional Laplacian have long been studied as interesting generalizations of the Brownian motion and its generator $\Delta$. Although jumping nature of stable processes of index $\alpha \in (0, 2)$ and non-local character of $\gen$ might suggest that their theory is quite different from the classical one, this is not the case. Recently much effort was made to obtain analogues of classical results concerning Dirichlet Laplacian, killed Brownian motion and harmonic functions in this ,,fractional'' setting. A large part of this work has already been done. It turned out that not only many statements remain true when the Brownian motion and $\Delta$ are replaced by the isotropic stable process $X(t)$ and $\gen$ with $\alpha \in (0, 2)$, but some can even be proved with much weaker regularity assumptions.

\smallskip

Consider, for instance, the so-called intrinsic ultracontractivity. It has been proved that the transition semigroup of the process $X(t)$ killed upon exiting an arbitrary open and bounded set $D$ is intrinsically ultracontractive, see~\cite{bib:k}, Theorem~1. This is false when the Brownian motion is considered, unless some regularity of the domain $D$ is assumed, e.g. that $D$ is a bounded and connected Lipschitz domain, cf.~\cite{bib:ds}, Theorems~9.1 and~9.3.

\smallskip

Among many other consequences, intrinsic ultracontractivity implies that the spectral gap determines the rate of convergence of transition density of the killed process to its ground state distribution. Hence, besides being interesting itself for analytic reasons, the estimates for the spectral gap have important probabilistic consequences. In this paper a lower bound is established for the spectral gap of the Dirichlet fractional Laplacian on an arbitrary bounded domain in $\Rd$ in terms of its diameter and the ground state eigenvalue. Clearly this result does not have a classical counterpart, since the spectral gap of the Dirichlet Laplacian on a (connected) domain may be arbitrarily small while its ground state eigenvalue as well as the diameter of the domain are bounded away from $0$ and $\infty$.

\medskip

Let us briefly recall some basic definitions. The \emph{isotropic stable process} of index $\alpha \in (0, 2)$ is a (pure jump) L\'evy process $X(t)$ with Fourier tranform $\ex_0 \exp(- i \left< z, X_t \right>) = \exp(-t \, |z|^\alpha)$; $\pr_x$ and $\ex_x$ stands for the probability distribution of the process starting at $x$ and the assosiated expected value. The isotropic stable process of index $2$ is the Brownian motion. The infinitesimal generator $\gen$ of $X(t)$ is given by
\[
  \gen f(x)
  =
  \frac{2^\alpha \, \Gamma(\frac{d+\alpha}{2})}
  {\pi^{\frac{d}{2}} \, |\Gamma(- \frac{\alpha}{2})|} \;
  \lim_{\varepsilon \searrow 0} \,
  \int\limits_{B(x, \varepsilon)^c}
  \frac{f(y) - f(x)}{|y - x|^{d + \alpha}} \, dy
  \, ,
\]
and is defined for an appropriate class of functions $f$. The symbol $\gen$ is short for $(-(-\Delta)^{\frac{\alpha}{2}})$, where the fractional power should be understood in terms of positive definite unbounded operators acting on $L^2(\Rd)$. The book of Dynkin~(\cite{bib:dy}) is an excellent introduction into the theory of Markov processes. See also~\cite{bib:s} for more specific information on L\'evy processes and~\cite{bib:bb} for information and references concerning isotropic stable processes.

\smallskip

Let $D$ be an open and bounded subset of $\Rd$. Consider the \emph{Dirichlet fractional Laplacian} on $D$, i.e. the negative Friedrichs extension of the operator mapping $f \in C_c^\infty(D)$ to $(- \gen f) \, \ind_D$. This is a generator of the isotropic stable process killed upon exiting $D$. It follows that there exists an orthonormal family of its eigenfunctions $\varphi_n \in L^2(D)$, corresponding to a nonincreasing sequence of negative eigenvalues $(-\lambda_n)$, see~\cite{bib:da}. All eigenfunctions are continuous on $D$ and bounded. The first eigenvalue is simple and, regarded as a function of the domain, monotone increasing in $D$. The corresponding eigenfunction $\varphi_1$, called the \emph{ground state eigenfunction}, may be chosen positive on $D$. \emph{Spectral gap} is the difference $(\lambda_2 - \lambda_1)$. As it was mentioned above, it describes how the transition density $p^D_t$ of the killed process converges to the ground state eigenfunction. More precisely,
\[
  \left| e^{\lambda_1 t} \, p^D_t(x, y) -
  \varphi_1(x) \, \varphi_1(y) \right|
  \le
  C \, e^{(\lambda_2 - \lambda_1) t} \, \varphi_1(x) \, \varphi_1(y)
\]
for $t \ge 1$, with some constant $C$ dependent on $\alpha$ and the domain $D$.

\begin{theorem}
\label{th:phi}
Let $\varphi_1$ be the ground state eigenfunction of the Dirichlet fractional Laplacian $\gen$ of any index $\alpha \in (0, 2)$ on an open and bounded set $D \sub \Rd$, $d \ge 1$. Denote by $(-\lambda_1)$ the corresponding eigenvalue of $\gen$. Then
\begin{equation}
\label{eq:phi}
  \sup \varphi_1
  \le
  c \, (\lambda_1)^{\frac{d}{2 \alpha}}
\end{equation}
with a constant $c$ dependent only on the dimension $d$ and index $\alpha$.
\end{theorem}

\noindent
From the proof of the theorem it follows that the constant is given by
\[
  c
  =
  \pi^{-\frac{d}{4}}
  \,
  \sqrt{2 d \, \Gamma(\mbox{$\frac{d}{2})$}}
  \left(\frac
    {4 \, \Gamma(\frac{d}{2})}
    {\alpha \, 2^\alpha \, \Gamma(\frac{d+\alpha}{2}) \,
    \Gamma(\frac{\alpha}{2})}
  \right)^{\frac{d}{2 \alpha}}
  \, .
\]
In particular, for the Cauchy process ($\alpha = 1$) in $\Rd$ the constant $c$ equals $2$ when $d = 1$ and $8 \pi^{-\frac{3}{2}} \le 1.44$ when $d = 2$. Notice that Theorem~\ref{th:phi} is also valid for the Dirichlet Laplacian on connected domains.

\begin{theorem}
\label{th:gap}
Denote by $(-\lambda_1)$ and $(-\lambda_2)$ two greatest eigenvalues of the Dirichlet fractional Laplacian $\gen$ of index $\alpha \in (0, 2)$ on an open and bounded set $D \sub \Rd$, $d \ge 1$. Then the spectral gap
\begin{equation}
\label{eq:gap}
  \lambda_2 - \lambda_1
  \ge
  \frac{\widetilde{c}}
  {(\lambda_1)^{\frac{d}{\alpha}} \, (\diam D)^{d + \alpha}}
\end{equation}
with a positive constant $\widetilde{c}$ dependent only on the dimension $d$ and index $\alpha$.
\end{theorem}

\noindent
As in the case of Theorem~\ref{th:phi}, the constant is given by an explicit formula
\[
  \widetilde{c}
  =
  \frac{2^\alpha \, \Gamma(\frac{d+\alpha}{2})}
  {c \, \pi^{\frac{d}{2}} \, |\Gamma(- \frac{\alpha}{2})|}
  \, .
\]
When $\alpha = 1$, the above reduces to $\widetilde{c} = \frac{1}{2 \pi} \ge 0.159$ if $d = 1$ and $\widetilde{c} = \frac{\sqrt{\pi}}{16} \ge 0.11$ if $d = 2$. It is worth noting that without further assumptions on the domain $D$, the degree of the estimate~\eqref{eq:gap} is optimal. This will be shown after proving Theorem~\ref{th:gap}. 

\begin{proposition}[\cite{bib:bk:1}, Corollary~2.2]
Suppose that $D$ is an open and bounded superdomain of a ball $B(x, r) \sub \Rd$, $d \ge 1$. The ground state eigenvalue $(-\lambda_1)$ of the Dirichlet fractional Laplacian $\gen$ on $D$, $\alpha \in (0, 2)$, satisfies
\begin{equation}
\label{eq:gs}
  \lambda_1
  \le
  \frac{\alpha \, (\alpha + \frac{d}{2}) \, \sqrt{\pi} \,
  \Gamma(\frac{\alpha}{2}) \, \Gamma(\alpha + \frac{d}{2})}
  {(a + d) \, \Gamma(\frac{1 + \alpha}{2}) \, \Gamma(\frac{d}{2})} \,
  r^{-\alpha}
  \, .
\end{equation}
\end{proposition}

\noindent
A proof of this easy result (valid also for Dirichlet Laplacian on connected domains) can by found at the end of the article. For an arbitrary open and bounded set $D$, Theorem~\ref{th:gap} and Proposition yield a numeric lower bound for the spectral gap of the Dirichlet fractional Laplacian on $D$. For instance, let $\alpha = 1$. For the interval $D = (-1, 1)$ the spectral gap exceeds $\frac{1}{3 \pi^2} \ge 0.033$. For the unit disk in $\R^2$ the lower bound equals $\frac{3}{256 \sqrt{\pi}} \ge 0.0066$ and for the square $D = (-1, 1)^2$ --- $\frac{3}{512 \sqrt{2 \pi}} \ge 0.0023$. These numbers are far away from exact values, e.g. it is known that the spectral gap $(\lambda_2 - \lambda_1)$ of $\Delta^{\frac{1}{2}}$ on $(-1, 1)$ is greater than $\lambda_1$, which in turn exceeds $1$, see~\cite{bib:bk:1}, Theorem~5.3 and its proof. From the other hand, the other two numeric estimates are better than any other known results, despite the emphasis of Theorem~\ref{th:gap} is on considering general open and bounded domains.

\smallskip

If $D$ is a convex planar domain symmetric with respect to both coordinate axes, then the spectral gap is greater than a constant multiple of $(\diam D)^{-2}$. Also precise asymptotics of $(\lambda_2 - \lambda_1)$ (with very rough constants) are known for rectangular domains, see~\cite{bib:dk}, Theorems~1.2 and~1.4.

\smallskip

From~\cite{bib:cs}, Example~5.1 it follows that the spectral gap for a convex and bounded domain $D$ is not less than $\frac{1}{2} \, (\widetilde{\lambda}_2)^{\frac{\alpha}{2}} - (\widetilde{\lambda}_1)^{\frac{\alpha}{2}}$, where $(-\widetilde{\lambda}_n)$ is the $n^{th}$ eigenvalue of the classical Dirichlet Laplacian on $D$. However, this lower bound is non-positive whenever $\alpha \le \log 4 / \log (\widetilde{\lambda}_2 / \widetilde{\lambda}_1)$. According to~\cite{bib:ab}, the ratio of $\widetilde{\lambda}_2$ to $\widetilde{\lambda}_1$ attains maximum equal to $(j_{\frac{d}{2}-1, 1})^2 / (j_{\frac{d}{2}, 1})^2$ when $D$ is a ball, $j_{p, k}$ being the $k^{th}$ positve zero of the Bessel function $J_p$. Therefore the estimate mentioned above is void in all dimensions $d \ge 4$, and even on the real line it is nontrivial only for $\alpha > 1$.

\smallskip

Theorem~\ref{th:gap} is hence new even for many convex, non-symmetric domains in $\R^2$, though it is certainly not optimal in this case. Furthermore it is also valid for arbitrary open and bounded sets, demonstrating the difference between classical and ,,fractional'' results mentioned in the first paragraph. The proof of Theorem~\ref{th:gap} avoids most of technical difficulties met in~\cite{bib:bkm},~\cite{bib:dk} and results concerning Dirichlet Laplacian, as it does not depend on geometric properties of the domain, such as convexity, conectedness or regularity of boundary.

\smallskip

There are many open problems concerning the shape of the ground state eigenfunction and the spectral gap. In~\cite{bib:bk:3},~\cite{bib:bkm} and~\cite{bib:dk} some very basic questions for convex domains are stated. For general open and bounded domains it would be interesting to describe the behaviour of level sets of the ground state eigenfunction $\varphi(x)$, as they play an important role in the proof of Theorem~\ref{th:phi} (see~\eqref{eq:lambda}). A better insight into the properties of $\varphi$ might also lead to more subtle versions of Theorem~\ref{th:gap}.


\section*{Proof of the results}

For an open and bounded $D$ let $\tau_D = \inf \{ t \ge 0 \, : \, X(t) \notin D \}$ be the \emph{first exit time}. Since $D$ is bounded, $\tau_D$ is finite a.s. Denote by $\omega^x_D$ the \emph{$\alpha$-harmonic measure} of $D$, i.e. the distribution of $X(\tau_D)$. Let $G_D$ be the Green's function of $D$ for the process $X(t)$, that is the kernel function of the Green's operator $G_D f(x) = \ex_x \int_0^{\tau_D} f(X_t) \, dt$ for continuous $f$. From analytic point of view, $G_D$ the inverse of the Dirichlet fractional Laplacian. Define $s_D(x) = \int G_D(x, y) \, dy$, so that $\gen s_D = -1$ on $D$. The probabilistic meaning of $s_D$ is the expected value of the exit time $\tau_D$, $s_D(x) = \ex_x \tau_D$. For the unit ball, $s_{B(0, 1)}(x) = 2^{1 - \alpha} \, \Gamma(\frac{d}{2}) \, (\alpha \, \Gamma(\frac{d+\alpha}{2}) \, \Gamma(\frac{\alpha}{2}))^{-1} \, (1 - |x|^2)^{\frac{\alpha}{2}}$, $x \in B(0, 1)$. From scaling properties of $\gen$ and $X_t$ it follows that $s_{r D}(x) = r^\alpha \, s_D(\frac{x}{r})$, $r D$ being the dilation of $D$.

\smallskip

The following isoperimetric-type result is a standard one; $|E|$ denotes the Lebesgue measure of $E$.

\begin{lemma}
Let $D \sub \Rd$ be open and bounded and $\alpha \in (0, 2)$. Suppose that $B = B(0, r)$ satisfies $|B| = |D|$. Then $s_D(x) \le s_B(0)$ for all $x \in D$.
\end{lemma}

\begin{proof}
With no loss of generality assume that $x = 0$. By Tonelli's theorem, $s_D(z) = \int_0^\infty \pr_z(\tau_D \ge t) \, dt$. Let $p_t(x)$ be the transition density of $X(t)$. Then
\[
\begin{split}
  \pr_0(\tau_D \ge t)
  & =
  \lim_{n \rightarrow \infty} \;
  \pr_0 \left(
    \forall \; 1 \le i < n \, t \; : \; X(\mbox{$\frac{i}{n}$}) \in D
  \right)
  \\ & =
  \lim_{n \rightarrow \infty} \;
  \idotsint \prod_{i = 1}^{\lfloor n \, t \rfloor} \ind_D(x_i) \,
  p_{\frac{1}{n}}(x_i - x_{i-1}) \, dx_1 \, ... \, dx_n
\end{split}
\]
with $x_0 = 0$. Since $p_t$ is a radial monotone decreasing function and the radial monotone decreasing rearrangement of $\ind_D$ equals $\ind_B$, an application of the Brascamp-Lieb-Luttinger inequality (see~\cite{bib:bll}, Theorem~3.4) yields that $\pr_0(\tau_D \ge t) \le \pr_0(\tau_B \ge t)$. Lemma follows.
\end{proof}

As a result of intrinsic ultracontractivity, $e^{\lambda_1 t} \, \pr_x(\tau_D \ge t)$ has a finite and positive limit $\varphi_1(x) \, \int \varphi_1(y) \, dy$, where $(-\lambda_1)$ and $\varphi_1$ are the first eigenvalue and the ground state eigenfunction of $\gen$ on $D$ respectively. Hence the inequality established in the proof of Lemma implies a similar isoperimetric-type result for $\lambda_1$: Among the sets of a fixed Lebesgue measure the ball has the greatest first eigenvalue $(-\lambda_1)$. The same argument yields analogous results in the classical case.

\begin{proof}[Proof of Theorem \ref{th:phi}]
For simplicity, denote $\varphi = \varphi_1$, $\lambda = \lambda_1$ and $M = \sup \varphi$. Let $U = \{ x \in D \, : \, \varphi(x) \ge M / 2 \}$. By continuity of $\varphi$ on $D$, the set $D \setminus U$ is open. Dynkin's formula (see~\cite{bib:dy}, Corollary to Theorem~5.1; also cf. \cite{bib:s}, Exercise~E.22 and \cite{bib:b}, p.~67) yields
\[
  \varphi(x)
  =
  \int\limits_{D \setminus U} \varphi(y) \, \omega^x_U(dy)
  +
  \lambda \int\limits_U G_U(x, y) \, \varphi(y) \, dy
  \, , \hspace{10mm} x \in U
  \, .
\]
Since $0 \le \varphi \le M / 2$ on $D \setminus U$ and $M / 2 \le \varphi \le M$ on $U$,
\[
  \lambda \, \frac{M}{2} \, \sup s_U
  \le
  \mbox{$\sup_U$} \, \varphi
  \le
  \frac{M}{2} + \lambda \, M \, \sup s_U
  \, .
\]
Recall that $\sup_U \varphi = M$. Hence
\begin{equation}
\label{eq:lambda}
  \frac{1}{2}
  \le
  \lambda \, \sup s_U
  \le
  2
  \, .
\end{equation}
The inequality $(M / 2)^2 \, |U| \le \int_U \, (\varphi(x))^2 \, dx \le 1$ yields that $M \le 2 \, |U|^{-\frac{1}{2}}$. Observe that $\sup s_U$ is estimated by $|U|$. Indeed, according to Lemma, among all open sets of a fixed Lebesgue measure the ball has the greatest supremum of the expected value of exit time. It follows that
\[
  |U|
  \ge
  C \, (\sup s_U)^{\frac{d}{\alpha}}
\]
with $C = (s_{B(0, 1)}(0))^{-\frac{d}{\alpha}} \, |B(0, 1)|$. Hence
\[
  M
  \le
  2 \, |U|^{- \frac{1}{2}}
  \le
  2 \, C^{- \frac{1}{2}} \, (\sup s_U)^{- \frac{d}{2 \alpha}}
  \, .
\]
The proof is completed by combining this with the lower bound of~\eqref{eq:lambda}.
\end{proof}

\begin{proof}[Proof of Theorem \ref{th:gap}]
Let $\varphi_1$ and $\varphi_2$ denote the ground state eigenfunction on $D$ and an eigenfunction corresponding to the second greatest eigenvalue $(-\lambda_2)$ respectively. Recently Dyda and Kulczycki proved that the spectral gap is given by the variational formula (see~\cite{bib:dk}, Theorem~1.1)
\[
\begin{split}
  \lambda_2 - \lambda_1
  & =
  C \, \inf_{f \in \mathcal{F}}
    \iint \frac{(f(x) - f(y))^2}{|x - y|^{d + \alpha}} \,
    \varphi_1(x) \, \varphi_1(y) \, dx \, dy
  \, ,
  \\
  \mathcal{F}
  & =
  \left\{
    f
    \; : \;
    \int (f(x) \, \varphi(x))^2 \, dx = 1
    \, , \;
    \int f(x) \, (\varphi(x))^2 \, dx = 0
  \right\}
  \, .
\end{split}
\]
Here $C = 2^{\alpha - 1} \, \pi^{-\frac{d}{2}} \, \Gamma(\frac{d+\alpha}{2}) \, |\Gamma(- \frac{\alpha}{2})|^{-1}$. In addition, the infimum is attained for $f = \frac{\varphi_2}{\varphi_1}$.

\smallskip

For simplicity, denote $R = \diam D$. Estimates $|x - y| \le R$ and $\varphi(z) \le M = \sup \varphi$, and the assertion following the variational formula yield
\[
  \lambda_2 - \lambda_1
  \ge
  \frac{C}{R^{d + \alpha} \, M^2}
  \inf_{f \in \mathcal{F}}
    \iint \left(
      \frac{\varphi_2(x)}{\varphi_1(x)} - \frac{\varphi_2(y)}{\varphi_1(y)}
    \right) ^2 \, (\varphi_1(x) \, \varphi_1(y))^2 \, dx \, dy
  \, .
\]
Since $||\varphi_1||_2 = ||\varphi_2||_2 = 1$ and $\left< \varphi_1, \varphi_2 \right> = 0$, the double integral equals $2$, so that
\[
  \lambda_2 - \lambda_1
  \ge
  \frac{2 \, C}{R^{d + \alpha} \, M^2}
  \, .
\]
An application of Theorem~\ref{th:phi} completes the proof.
\end{proof}

Consider the set $D = B(-a, 1) \cup B(a, 1)$, where $a = (r, 0, ..., 0)$ and $r > 2$. Since $D$ is symmetric with respect to the reflection changing the sign of the first coordinate, so is $\varphi_1$. By the variational formula with $f = \ind_{B(a, 1)} - \ind_{B(-a, 1)}$ and by Theorem~\ref{th:phi},
\[
  \lambda_2 - \lambda_1
  \le
  \int\limits_{B(-a, 1)} \int\limits_{B(a, 1)}
    \frac{8 \, C}{r^{d + \alpha}} \, \varphi(x) \, \varphi(y) \, dx \, dy
  \le
  C' \, (\lambda_1)^{\frac{d}{\alpha}} \, r^{-d - \alpha}
  \, .
\]
The constant $C'$ depends only on $d$ and $\alpha$, and $\lambda_1$ does not exceed the ground state eigenvalue of $(-\gen)$ restricted to the unit ball (recall that $\lambda_1$ is a decreasing function of the domain). This implies that the degree of the estimate~\eqref{eq:gap} cannot be improved, justifying the note following Theorem~\ref{th:gap}.

\begin{proof}[Proof of Proposition]
Observe that $\left< -\gen f, f \right> \ge \lambda_1 ||f||_2^2$ for any $f \in L^2(D)$. Let $f = s_{B(x, r)}$, so that $\gen f = -1$ on $B(x, r)$. The explicit formula for $s_{B(x, r)}$ and integration in spherical coordinates yield
\[
  \lambda_1 \,
  \int_0^r t^{d-1} \, (r^2 - t^2)^\alpha \, dt
  \le
  \frac{2^{\alpha - 1} \, \alpha \, \Gamma(\frac{d+\alpha}{2}) \,
  \Gamma(\frac{a}{2})}{\Gamma(\frac{d}{2})} \,
  \int_0^r t^{d-1} \, (r^2 - t^2)^{\frac{\alpha}{2}} \, dt
  \, .
\]
Estimate~\eqref{eq:gs} follows by a change of variables $s = \frac{t^2}{r}$ and beta integration.
\end{proof}



\end{document}